\documentclass[10pt, reqno]{amsart}
\usepackage[arrow,matrix,curve]{xy}

\newtheorem{theorem}{Theorem}[section]

\newtheorem{lemma}[theorem]{Lemma}
\newtheorem{corollary}[theorem]{Corollary}
\newtheorem{conjecture}[theorem]{Conjecture}

\theoremstyle{definition}

\newcommand{\R}{\mathbb R} 

\newcommand{\RP}{\mathbb RP}

\numberwithin{equation}{section}
\numberwithin{figure}{section}

\begin{document}

\author[G.~Chambers]{Gregory R. Chambers$^1$} \address{ Department of
Mathematics, University of Chicago, Chicago, IL, 60637-1505 USA} \email{chambers@math.uchicago.edu}
\thanks{$^1$Supported by an NSERC Postdoctoral Fellowship}

\author[C.~Croke]{Christopher Croke} \address{ Department of
Mathematics, University of Pennsylvania, Philadelphia, PA
19104-6395 USA} \email{ccroke@math.upenn.edu}

\author[Y.~Liokumovich]{Yevgeny Liokumovich} \address{ Department of
Mathematics, MIT, Cambridge, MA 02139-4307} \email{ylio@mit.edu}

\author[H.~Wen]{Haomin Wen$^2$} \address{ Foresee Fund, Shanghai, 200122 China} \email{wenhaomin@foreseefund.com}
\thanks{$^2$Supported by the Max-Planck
Inst. Bonn 2014-2015}


\title[hemispheres] {Area of convex disks}


\begin{abstract}
This paper considers metric balls $B(p,R)$ in two dimensional Riemannian manifolds when $R$ is less than half the convexity radius.  We prove that
$Area(B(p,R)) \geq \frac 8 \pi R^2$.  This inequality has long been conjectured for $R$ less than half the injectivity radius.  This result also yields the upper bound $\mu_2(B(p,R)) \leq 2(\frac{\pi}{2R})^2$ on the first nonzero Neumann eigenvalue $\mu_2$ of the Laplacian in terms only of the radius.  This has also been conjectured for $R$ up to half the injectivity radius.

\end{abstract}

\maketitle

\section{Introduction}

In this short note we will consider balls $B(p,R)$ in complete two dimensional Riemannian manifolds $M$ for sufficiently small $R$.  We will let $inj(M)$ be the injectivity radius of $M$ and $conv(M)$ the convexity radius of $M$.  By definition $conv(M)$ is the smallest number such that for all $R < conv(M)$,  $B(p,R)$ is strictly convex (i.e for every $x,y \in B(p,r)$ there is a unique minimizing geodesic from $x$ to $y$ and it lies in $B(p,r)$).  Note that this implies that the boundary curve is convex.  Also note that $conv(M)\leq \frac 1 2 inj(M)$.

There is a long standing conjecture, in all dimensions $n$, that
hemispheres have the smallest volume among balls $B(p,R)$ of a fixed
dimension $n$ and radius $R \leq \frac 1 2 inj(M)$ . By a hemisphere we will mean a ball $B(x,R)$
of radius $R$ in the sphere with constant curvature $(\frac  \pi{2R})^2$.
For example, when $R=\frac \pi 2$ this is isometric to a
hemisphere of the unit sphere.

\begin{conjecture}
If $R \leq \frac 1 2 inj(M)$  then
$$ Vol(B(p,R))\geq \frac {\alpha(n)} 2\big( \frac {2R} \pi \big)^n$$
where $\alpha(n)$ represents the volume of the unit $n$-sphere.
Further equality holds if and only if $B(p,R)$ is isometric to a
hemisphere of (intrinsic) radius $R$.
\end{conjecture}

Although in all dimensions there are known (nonsharp) constants
$C(n)$ such that $Vol(B(p,R))\geq C(n)(\frac {2R} \pi)^n$ (see
\cite{Be} for $n=2,3$ and \cite{Cr} for all $n$) even the two
dimensional case of the conjecture is open:

\begin{conjecture}
\label{sareaball} For a two dimensional surface $M$.  If $R \leq \frac 1 2 inj(M)$ then
$$Area(B(p,R)) \geq \frac 8 \pi R^2.$$
Further equality holds if and only if $B(x,R)$ is isometric to a
hemisphere of (intrinsic) radius $R$.
\end{conjecture}

The best known result of this type (see \cite{Cr1}) is that $Area(B(p,R) \geq \frac {8-\pi}{2} R^2$.

Although we cannot solve this conjecture for the full range of $R$ we show it is true for $R\leq \frac 1 2 conv(M)$.  (Hence in our case $R$ is always $\leq \frac 1 4 inj(M)$.)

\begin{theorem}
\label{Main}
For a two dimensional surface $M$.  If $R \leq \frac 1 2 conv(M)$ then
$$Area(B(p,R)) \geq \frac 8 \pi R^2.$$
\end{theorem}

\bigskip
\bigskip

For a ball $B(p,R)$ let $0<\lambda_1\leq\lambda_2\leq
\lambda _3,...$ be the spectrum of the Laplace operator with
Dirichlet boundary conditions and $0=\mu_1<\mu_2\leq \mu_3\leq...$
be the spectrum for Neumann boundary conditions.  The unit
hemisphere (i.e. a ball of intrinsic radius $\frac \pi 2$ in the
unit sphere) has $\lambda_1=\mu_2=\mu_3=2$. Thus a hemisphere of
intrinsic radius $R$ has $\lambda_1=\mu_2=\mu_3=2(\frac \pi
{2R})^2$ where the corresponding eigenfunctions are the coordinate
functions (from the embedding in $R^3$).

In \cite{Cr1} it is shown how the estimate in the theorem along with results in \cite{Cr} and \cite{He} yield the following.

\begin{corollary}
\label{seigen} For $B(p,R)$ a metric ball on a surface with $R\leq \frac 1 2 conv(M)$ then
$$\mu_2\leq 2(\frac \pi {2R})^2\leq \lambda_1.$$
\end{corollary}

The corresponding result for $R\leq \frac 1 2 inj(M)$ is conjectured to be true in \cite{Cr1}.

\bigskip
\bigskip

Another open problem that has arisen in the context of these questions is: ``must closed geodesic triangles in $B(p,R)$ that wind around $p$ and have vertices on the boundary have length $\geq 4R$?''.  We show below (Corollary \ref{polygon}) that this is true for $R\leq \frac 1 2 conv(M)$.  The question is still open for $R\leq \frac 1 2 inj(M)$.

\section{The Proof}

The tool we use from the convexity property is the following.

\begin{lemma}
Let $\gamma(t)$ be a unit speed geodesic contained in a ball $B(p,R)$ of radius $R\leq conv(M)$ in a Riemannian surface $M$.  Then the function
$f(t)=d(p,\gamma(t))$ is convex.
\end{lemma}

{\bf Proof:} If $\gamma(t_0)=p$ then $f(t)=|t-t_0|$ is convex, so we can assume that $\gamma$ does not pass through $p$.  On $B(p,R) - \{p\}$ let $\vec r = \frac {\partial}{\partial r}$ be the radial vector field and let $\vec \theta = \frac {\partial}{\partial \theta}$.  Let $g$ be the function on $B(p,R) - \{p\}$ defined by
$$g= - \langle \nabla_{\vec \theta}{\vec \theta}, \vec r \rangle =  \langle \nabla_{\vec \theta}{\vec r}, \vec \theta \rangle .$$
The convexity of the boundary of $B(p,r)$ for $r\leq R$ says that $g\geq 0$.  Let $\gamma'(t)=a(t)\vec r+b(t)\vec \theta$.  Since $f'(t)=\langle \gamma'(t),\vec r \rangle$ we can compute;
$$f''(t)= \langle \gamma'(t),\nabla_{\gamma'(t)}\vec r \rangle=\langle a(t)\vec r+b(t)\vec \theta, b(t)\nabla_{\vec \theta}{\vec r} \rangle =g(\gamma(t))\cdot b(t)^2 \geq 0.$$
\qed

From this we get the following inequality.

\begin{corollary}
\label{symmdist}
Let $B(p,R)$ be a ball of radius $R\leq \frac 1 2 conv(M)$ in a Riemannian surface $M$.  Let $q$ be a boundary point of $B(p,R)$ and $\gamma(t)$ be a unit speed geodesic with $\gamma(0)=p$.  Then
$$d(q,\gamma(t))+d(q,\gamma(-t))\geq 2R.$$
\end{corollary}

{\bf Proof:}  This follows pretty directly from the lemma.  Since $f(t)=d(q,\gamma(t))$ is convex, $f(t)+f(-t)\geq 2f(0)=2R$ yielding the result.

\qed

\begin{corollary}
\label{polygon}
Let $R\leq \frac 1 2 conv(M)$ and $T\subset B(p,R)$ be a geodesic triangle with vertices on $\partial B(p,R)$.  If $p \notin T$ we assume that $T$ winds around $p$ (i.e. has nonzero winding number in $B(p,R)$).  Then $L(T)\geq 4R$.
\end{corollary}

{\bf Proof:} Since the result is clear if $T$ passes through $p$ we can assume it does not.  Diameters (i.e. geodesics passing through $p$) intersect $T$ at exactly two points with $p$ between them on the diameter.  A simple intermediate value theorem argument shows that there is a diameter $\gamma$ such that the two points are equidistant from $p$ (i.e. $\gamma(t_0)$ and $\gamma(-t_0)$ when $\gamma(0)=p$).  Thus $\gamma$ breaks $T$ into two parts each of which has at least one vertex.  It follows from Corollary \ref{symmdist} and the triangle inequality that
each part has length at least $2R$, yielding the result.

\qed

Note we could consider simple closed geodesic polygons $P$ with vertices on the boundary and that wind around $p$. The above along with straightforward triangle inequalities show that $L(P)\geq 4R$.

\bigskip

{\bf Proof of theorem \ref{Main}:}  Let $\gamma$ be a diameter, i.e. a geodesic through the center $p$.  Then $\gamma$ divides $B(p,R)$ into two halves $B_+$ and $B_-$.  We will prove the result by showing (with the same argument) that each has Area $\geq \frac 4 \pi R^2$.  We will do this by constructing a metric on $\RP^2$ and applying Pu's theorem \cite{Pu}.  First construct a ($C^1$-smooth) metric $B_{++}$ on the disc by gluing two copies, $B^1_+$ and $B^2_+$, of $B_+$ together along the two copies, $\gamma^1$ and $\gamma^2$, of the $\gamma$ portion of the boundaries by identifying $\gamma^1(t)$ with $\gamma^2(-t)$.  If $i:B_+\to B_+$ is the geodesic inversion through $p$ then by construction it is an isometry.  The area of $B_{++}$ is twice the area of $B_+$.

We claim that $d(q,i(q))=2R$ for any $q\in \partial B_{++}$.  Say $q\in B^1_+$ and $\tau:[0,1] \to B_{++}$ any path from $q$ to $i(q)$.  There must be a $t_0\in [0,1]$ such that $\tau(t_0)=y\in \gamma$.  Let $s_0$ be such that $y=\gamma^1(s_0)=\gamma^2(-s_0)$.  Now the shortest path from $y$ to $q$ is the geodesic $\sigma_1$ in $B^1_+$ between $y$ and $q$.  (Note that if a path from $y$ to $q$ ever left $B^1_+$ there is a shorter path replacing the parts in $B^2_+$ with segments of $\gamma$.)  Similarly the shortest path from $y$ to $i(q)$ is the geodesic $\sigma_2$ in $B^2_+$.  Thus, by the Corollary \ref{symmdist}, $L(\tau)\geq L(\sigma_1)+L(\sigma_2)=L(\sigma_1)+L(i(\sigma_2))\geq 2R$.  Thus $d(q,i(q))\geq 2R$.  The claim follows since the geodesic through $p$ from $q$ to $i(q)$ has length $2R$.

We now consider the metric (only $C^0$ along $\partial B_{++}$) on $\RP^2$ obtained by identifying $q\in \partial B_{++}$ with $i(q)$.  We claim that the systole (the length of the shortest non-contractible closed curve) is equal to $2R$.  Certainly the geodesics in $B_{++}$ through $p$ from $q$ to $i(q)$ become non-contractible and hence the systole is $\leq 2R$.  For any closed curve there is a nearby curve of nearly the same length that intersects $\partial B_{++}$ transversely in finitely many points.  If $\tau:[0,1] \rightarrow \R P^2$ is a closed non-contractible loop, then it must intersect $\partial B_{++}$ at least once.  We can assume that $\tau(0)=\tau(1)\in B_{++}$ and $\tau$ intersects $\partial B_{++}$ transversely finitely often (i.e. there are $0=t_0<t_1<...t_n=1$ such that $\tau(t_i)\in \partial B_{++}$).  Note that $\tau$ need not be a continuous curve when thought of in $B_{++}$.  However, by replacing every other segment $\tau ([t_{2i-1},t_{2i}])$ with $i(\tau ([t_{2i-1},t_{2i}]))$, we create a new curve $\bar \tau$ of the same length as $\tau$ that is homotopic in $\RP^2$ to $\tau$ and is continuous as a curve in $B_{++}$.  Further, since $\tau$ was not contractible, $\bar \tau(1)=i(q)$.  The previous claim, however, says that $2R\leq L(\bar \tau)=L(\tau)$, and so the systole is $2R$.

We now apply Pu's theorem \cite{Pu} to conclude that $\frac 8 \pi R^2\leq Area(\R P^2)=Area(B_{++})=2Area(B_+)$ and the theorem follows.  Although Pu's Theorem assumes smooth metrics, our metrics can be approximated by smooth metrics with nearly the same volume and systole so the inequality applies to our metrics as well.  The same approach proves that $Area(B_{-}) \geq \frac{4}{\pi}R^2$ as well, completing the proof.

\qed

\end{document}